\magnification=1200
\input amstex
\documentstyle{amsppt}
\addto\tenpoint{\normalbaselineskip13pt\normalbaselines}
\pageheight{44pc}
\pagewidth{30pc}

\font\got=eufm10

\message{version of February 9, 1999}


\newbox\Text
\newcount\width
\newcount\height
\newdimen\dimh
\newdimen\dimv
\newcount\dep
\newdimen\shifth
\newdimen\shiftv
\newcount\ellmarg

\newcount\cx
\newcount\cy

\newbox\O
\newdimen\dimhO
\newdimen\dimvO
\newcount\hO
\newcount\vO

\newbox\FText
\newdimen\marg
\newdimen\rs
\newdimen\scrspc


\def\encirc#1{
    \ellmarg=30         
    \rs=1pt             
    \scrspc=-0.5pt      
    \setbox\O=\hbox{$\bigcirc$}
    \dimhO=\wd\O
    \dimvO=\ht\O
    \hO=\dimhO
    \vO=\dimvO
    \divide\hO by 9472
    \divide\vO by 9472
    \advance\hO by \ellmarg
    \advance\vO by \ellmarg
    \advance\hO by -19
    \advance\vO by -10
    \setbox\Text=\hbox{$\scriptstyle{#1}\scriptspace=\scrspc$}
    \dimh=\wd\Text
    \dimv=\ht\Text
    \dep=\dp\Text
    \divide\dep by 9472
    \width=\dimh
    \divide\width by 9472
    \cx=-\width
    \advance\width by \ellmarg
    \height=\dimv
    \divide\height by 9472
    \cy=-\height
    \advance\cy by \dep
    \advance\height by \dep
    \advance\height by \ellmarg
%
%
    \advance\hO by -\width
    \advance\vO by -\height
    \ifnum\hO>0
        \ifnum\vO>0
            \putcircle
        \else
            \putellipse{\width}{\height}{\cx}{\cy}
        \fi
    \else
        \putellipse{\width}{\height}{\cx}{\cy}
    \fi
}

\def\putellipse#1#2#3#4{
\mskip8mu \raise\rs\hbox{%
\copy\Text%
\special{pn 5}%
\special{ar \the#3 \the#4 \the#1 \the#2 0 6.28}%
}\mskip8mu
}


\def\putcircle{
\shifth=\dimhO
\advance\shifth by \dimh
\divide\shifth by 2
\mskip3mu\hbox{\copy\O\kern-\the\shifth\raise\rs\hbox{\copy\Text}}\mskip7mu
}



\def\1{%
\hbox{1\kern-3.5pt 1\kern-4pt%
\raise 5.7pt\hbox{\vrule height1pt width1.5pt depth0pt}}
}


\def\ensquare#1{
    \rs=1pt 
    \marg=1pt 
    \scrspc=-0.5pt      
    \setbox\Text=\hbox{$\scriptstyle{#1}\scriptspace=\scrspc$}
    \setbox\FText=\hbox{\vrule%
        \vbox{\hrule\kern\marg%
              \hbox{\kern\marg\copy\Text\kern\marg}%
              \kern\marg\hrule}%
        \vrule}
    \shiftv=\dp\Text
    \advance\shiftv by \marg
    \advance\shiftv by -\rs
    \mskip4mu\lower\shiftv\hbox{\unhbox\FText}\mskip5mu
}


\TagsOnRight

\def\endthm{\endproclaim}
\define\df#1 {\definition{Definition #1} }
\def\enddf{\enddefinition}
\define\ex#1{\example{Example #1}}
\def\endex{\endexample}
\define\endlemma{\endproclaim}
\def\endprop{\endproclaim}
\define\rem#1 {\remark{Remark \rom{#1}}}
\def\endrem {\endremark}
\define\pf {\demo{Proof}}
\def\endpf{\qed\enddemo}
\def\endcor{\endproclaim}
\def\ms {\medskip}

\def\gk {\text {GKdim\,}}
\def\coef {\text {Coeff\,}}
\def\span {\text {Span}}
\def\rk {\text {rk\,}}
\def\ad {\text {ad\,}}

\def\cur {\text {Cur}}
\def\Mat {\Cal M\text{at}}

\def\s#1 {_{(#1)}}
\def\c#1#2#3{{#1}\encirc{#2}{#3}}
\def\d {\partial}
\def\binom#1#2{{#1 \choose #2}}
\def\C {{\Bbb C}}
\def\Z {\Bbb Z}
\def\g {\hbox{\got g}}

\topmatter
\title
On Associative Conformal Algebras of Linear Growth
\endtitle
\author Alexander Retakh\endauthor
\address Department of Mathematics, Yale University, New Haven, CT, 06520
\endaddress
\email retakh\@math.yale.edu\endemail
\abstract
We introduce the notions of conformal identity and unital associative 
conformal algebras and classify finitely generated simple unital associative 
conformal algebras of linear growth.  These are precisely the complete 
algebras of conformal endomorphisms of finite modules.
\endabstract

\endtopmatter

\document

\head 0. Introduction.\endhead

The subject of conformal algebras is a relatively recent development in
the theory of vertex algebras \cite{K1}.  The relation of Lie conformal
algebras to vertex algebras is similar to that of Lie algebras and their
universal enveloping algebras.

Semisimple Lie conformal algebras of finite type were classified in
\cite{DK} and semisimple associative algebras of finite type in \cite{K2}.  
Associative conformal algebras appear in conformal representation theory.  
The complete algebras of conformal endomorphisms of finite modules
$\text{CEnd}_N$ (called in this paper the conformal Weyl algebras) are a
particular example(see \cite{K1, 2.9}).  Notice that these are algebras
are not finite as modules over $\C[\d]$.

These paper is concerned with associative algebras of finite growth (but
not necessarily of finite type). We require the conformal algebra to be
unital, that is, to possess an element that acts as a left identity with
respect to the $0$-th multiplication and whose locality degree with itself
is $1$.  In particular, this means that its $0$-th coefficient is the
(left) identity of the coefficient algebra.  One can then use the
classification of associative algebras of linear growth obtained in
\cite{SSW} to classify a class of unital conformal algebras.

\proclaim{Main Theorem} Let $C$ be a finitely generated simple unital
associative conformal algebra. If $C$ has growth $1$, then it is a
conformal Weyl algebra $\text{CEnd}_N$.\endproclaim

Notice that there is a marked difference from the case of growth $0$ where
simple associative conformal algebras are just the current algebras over
the matrix algebras (\cite{K2, 4.4}, see also Remark 4.3 below).

The paper is organized as follows.  The first chapter is devoted to
preliminary material on conformal algebras where, for the most part, we
loosely follow the treatment of \cite{K2}.  The second chapter discusses
the concept of Gelfand-Kirillov dimension for conformal algerbas and
relates it to the dimension of the coefficient algebra.  The third chapter
is devoted to the discussion of unital conformal algebras.  Here the main
result (Proposition 3.5) relates the presence of conformal identity to the
structure of the coefficient algebra. The fourth chapter contains the
proof of the Main Theorem in Theorem 4.6.  We also classify unital
associative algebras of growth $0$ in Theorem 4.2 (the proof is different
from \cite{K2} which contains the general case).

I am extremely grateful to Efim Zelmanov for introducing me to the subject
of conformal algebras and guiding me through all stages of this research
and to Michael Roitman for helpful discussions and his circle-drawing \TeX
package.

\head 1. Preliminaries on Conformal Algebras.\endhead

We start with a motivation for the concept of the conformal algebra.
Let $R$ be a (Lie or associative) algebra.  One can consider formal
distributions over $R$, that is elements of $R[[z,z^{-1}]]$.  These
appear, for instance, in the theory of  operator product expansions in
vertex algebras.
Obviously, it is impossible to multiply two formal distributions; however,
one can introduce an infinite number of bilinear operations that act
as a ``replacement'' for multiplication.  Let $a(z), b(z)$ be two
formal distributions.  For an integer $n\geq 0$ define another formal 
distribution 
$$\c {a(z)}n{b(z)}={\hbox {\rm Res}}|_{z=0}\,a(z)b(w)(z-w)^n\tag 1.1$$ 
called the $n$-th 
product of $a(z)$ and $b(z)$.  (Those familiar with vertex algebras will
observe that we do not consider the $-1$st, i.e. the Wick, product here
as its definition requires a representation of $R$.)  If one writes
$a(z)$ explicitly as $\sum a(n)z^{-n-1}$, it follows that 
$$(\c {a(z)}n{b(z)})(k)=
\sum_{j\in\Z_+} (-1)^j \binom nj a(n-j)b(k+j).\tag 1.2$$

In the theory of vertex algebras one wants to consider only the mutually
local series, i.e. $a(z)$ and $b(z)$ such that $\c {a(z)}n{b(z)}$ and $\c
{b(z)}n{a(z)}$ are zero for $n>>0$. We simply say that $a(z)$ and $b(z)$
are local if $\c {a(z)}n{b(z)}=0$ for $n>N(a,b)$ and call minimal such
$N(a,b)$ the degree of locality of $a(z)$ and $b(z)$.  In terms of
coefficients this becomes:
$$\sum_{j=0}^n (-1)^j \binom nj a(n-j)b(k+j)=0,\ \ \ n>>0, k\in\Z.\tag
1.3$$ It is possible to rewrite this statement as follows:

\proclaim{Lemma 1.1} \cite{K1, 2.3} The series $a(z)$ and $b(z)$ are
local iff
$a(z)b(w)(z-w)^n=0$ for $n>N(a,b)$.  In terms of 
coefficients:
$$\sum_{j=0}^n (-1)^j\binom nj a(l-j)b(m+j)=0,\ \ \
n>N(a,b), l,m\in\Z.\tag 1.4$$\endlemma

Apart from locality one also wants to take into account the action of
$\d=\d/\d z$.  These consideration suggested the following definition,
first stated in \cite{K1}.

\df{1.2} A conformal algebra $C$ is a $\C[\d]$-module endowed with bilinear
products $\encirc{n},\, n\in \Z_{\geq 0}$, that satisfy the following
axioms for any $a,b\in C$:

(C1) (locality) $\c anb=0$ for $n>N(a,b)$;

(C2) (Leibnitz rule) $\d(\c anb)=\c {(\d a)}nb+\c an{(\d b)}$;

(C3) $\c {(\d a)}nb=-n\c a{n-1}b$.
\enddf

Clearly algebras of formal distributions closed with respect to the action
of $\d$ satisfy these axioms. 

The bilinear product $\encirc{n}$ is usually called the multiplication of
order $n$.
\ms

In general, $N(a,b)\neq N(b,a)$.  Moreover, even if $a,b$ and $c$ are
pairwise local, $\c anb$ and $c$ need not be such.  One can, however, 
establish some 
correlations between different degrees of locality.  For instance, it
follows from (C3) that $N(\d a,b)=N(a,b)+1$ and from (C2) and (C3) that
$N(a, \d b)=N(a,b)+1$. It follows as well that
$\c {(\d a)}0b=0$ for all $a,b$.  Yet, one can hardly make deeper
statements about conformal algebras without restricting attention to
less general cases.

Two such cases are algebras of formal distributions over Lie and
associative algebras closed with respect to the derivation $\d$.  Here
taking products does not destroy locality:

\proclaim{Lemma 1.3} {\bf (Dong's lemma)} \cite{Li, K1} Let $a,b$ and $c$
be pairwise mutually local formal distributions over either a Lie or
associative algebra. Then for any $n\geq 0$ $\c anb$ and $c$ are again
pairwise mutually local. \endlemma

We remarked above that $\d$ preserves locality (although not its degree)
as well.  Thus for such algebras of formal distributions it is enough
to check locality of the generators.  Examples follow:

\ex{1.4} 1. {\it Current algebras}: let $\g$ be a Lie algebra and
$\g[t,t^{-1}]$
its affine extension; $[gt^m,ht^m]=[g,h]t^{m+n}$.  The algebra of formal
distributions generated by $g(t)=\sum gt^nz^{-n-1}$ and their derivatives
is a conformal algebra.  The non-zero products for generators are $$\c
{g(t)}0{h(t)}=[g,h](t).$$  This algebra is denoted $\cur(\g)$.
\ms\par 
2. {\it (conformal) Virasoro algebra}: let $t, t^{-1}, \d_t$ be the
generators of the Lie algebra of differential operators on $S^1$ (i.e. a
centerless Virasoro algebra).  Consider the conformal algebra generated by
the formal distribution $L=-\sum \d_tt^nz^{-n-1}$ and its derivatives.  
(This is a non-standard way of writing the generator of the conformal
Virasoro algebra and, for that matter, of the centerless Virasoro algebra.  
We put powers of $t$ on the left since later it will become the generating
indeterminate of a polynomial extension; see below.) The non-zero products
for the generator are $$\c L0L=\d L \text{\ and\ } \c L1L=2L.$$ \par
3. {\it Associative current algebras}: let $A$ be an associative algebra
and $A[t,t^{-1}]$ its affine extension ($t$ commutes with $A$).  The
algebra of formal distributions generated by $f_a=\sum at^nz^{-n-1}$ and
its derivatives is a conformal algebra.  The non-zero products for the
generators are $$\c {f_a}0{f_b}=f_{ab}.$$ This algebra is denoted
$\cur(A)$.
\ms\par
4. {\it (conformal) Weyl algebra} ${\hbox {\got W}}$: let $t, \d_t$ be the
generators of the Weyl algebra $W$ with the relation $\d_tt-t\d_t=1$.
Consider its localization at $t$ denoted here by $W_t$. Formal
distributions $e=\sum t^nz^{-n-1}$ and $L=\sum \d_tt^nz^{-n-1}$ with their
derivatives generate a conformal algebra. The non-zero products for the
generators are $$\eqalign{&\c e0e=e, \c e0L=\c L0e=L,\cr &\c e1L=\c
L1e=-e,\cr &\c L0L=\d_t L, \c L1L=-L.}$$ The algebra of formal
distributions over $\Mat_n(W_t)$, which is a generalization of the above
example, will also be called a (conformal) Weyl algebra. \endex

Examples 1 and 2 were studied extensively in \cite{DK}.  Obviously, the
third example mimics example 1.  The Weyl algebra and the Virasoro algebra
are also closely related.  Since any associative conformal algebra can be
turned into a Lie conformal algebra (\cite{K1, 2.10a}), the Weyl algebra
will contain a subalgebra isomorphic to the Virasoro algebra (and $e$ will
become a central element). Again, let us remark that the generator of the
Virasoro algebra is usually written as $L=-\sum t^n\d_tz^{-n-1}$ (the
results of conformal multiplications are the same here).  In a similar
manner, one can choose the formal distribution $\sum t^n\d_tz^{-n-1}=L+\d
e$ as a generator of the Weil algebra instead of $L$.

Another approach to the Weyl algebra (or rather the matrix algebra over
the Weyl algebra) is to consider the complete algebra of conformal
endomorphisms of a finite $\C[\d]$-module.  For a discussion of conformal
representation theory see \cite{DK} and \cite{K2}.

It is possible to generalize example 4 by considering a conformal algerba
of series in $A[t, t^{-1};  \delta][[z,z^{-1}]]$ where $A$ is an
associative algebra, $t$ an indeterminate, and $\delta$ a locally
nilpotent derivation of $A$.  Such algebras are called {\it differential
algebras}. We postpone their discussion until Chapter 3.
\ms

One obviously wants to formalize the above discription of Lie and
associative algebras (and, in particular, Dong's lemma 1.3) of formal
distributions to the general case of conformal algebra.  

Remark, first, that the $0$-th product of a Lie or an associative algebra
of formal distributions behaves just as a regular Lie or associative
product.  It is certainly possible to forget about products of positive
order and consider the algebra $(C,\encirc{0})$.  However, one wants to
take into account the full structure of the conformal algebra, not just
its $0$-th product. 

The following construction was first partially introduced in
\cite{Bo} for vertex algebras and generalized in \cite{K1} for conformal
algebras. 

Notice first that by axioms (C2) and (C3), $\d C$ is an ideal in
$(C,\encirc{0})$.
Consider now a conformal algebra $\tilde C=C[t,t^{-1}]$.  Define 
the $n$-th product as
$$\c {at^l}n{bt^m}=\sum_{j\in\Z_+} \binom lj(\c a{n+j}b)t^{l+m-j}.\tag
1.5$$ and put $\tilde \d=\d+ \d/\d t$.
This makes $\tilde C$ into a conformal algebra.

\df{1.5} (\cite {K1}) The coefficient algebra $\coef C$ of conformal
algebra $C$ is the quotient algebra $(\tilde C,\encirc{0})/\tilde\d\tilde
C$.\enddf

Denote the map $\tilde C \to \coef C$ by $\phi$. Consider the algebra of
formal distributions over $\coef C$ consisting of series $\sum
\phi(at^n)z^{-n-1}$ where $a\in C$.  This algebra is isomorphic to $C$. 
It obviosuly follows that every conformal algebra can be written as an
algebra of formal distributions over its coefficient algebra.

In the remaining chapters we will often write $\phi(at^n)$ as simply
$a(n)$, i.e. we will not distinguish between a conformal algebra and the
corresponding algebra of formal distributions.  The following formula will
be quite useful: $$(\d a)(n)=-na(n-1), \ \ \ n\in\Z.\tag 1.6$$

We say that $C$ is conformal Lie or associative if $\coef C$ is respectively
Lie or associative (this definition can be applied to any variety of
algebras as well).  In particular a conformal algebra is Lie conformal if
$$\eqalign{\c fng&=\sum_{j\in\Z_+}{(-1)^{j+n}\over j!}\d(\c g{n+j}f),\cr
\c fm{(\c gnh)}&=\sum_{j\in\Z_+}\binom mj (\c fjg\c ){m+n-j}h+\c
gn{(\c fmh)}}\tag 1.7$$
and associative conformal if
$$\c fm{(\c gnh)}=\sum_{j\in\Z_+}\binom mj(\c fjg\c ){m+n-j}h.\tag 1.8$$
The latter condition is equivalent to
$$(\c fmg\c )nh=\sum_{j\in\Z_+}(-1)^j\binom mj \c f{m-j}{(\c
g{n+j}h)}.\tag 1.9$$

Many standard definitions from the ``classical'' structure theory carry
over to the conformal case in an obvious manner.  Thus, a conformal ideal
is a conformal subalgebra closed under all right and left multiplications
of any order, a simple conformal algebra is an algebra that contains no
ideals, and so on. 

\ms

Below all conformal algebras are assumed to be associative conformal
unless otherwise stated.

\ms

\head 2. Gelfand-Kirillov Dimension of Conformal Algebras.\endhead

The Gelfand-Kirillov dimension of a finitely generated algebra (of any
variety) $A$ is defined as $$\gk A=\limsup_{r\to\infty}{\log {\dim}
(V^1+V^2+\dots+V^r)
\over \log r},$$ where $V$ is a generating subspace of $A$.  This
definition easily carries over to the conformal case.  

Let $C$ be a finitely generated conformal algebra (over any variety) with
generators $f_1,\dots,f_n$.  Define $C_r$ to be a $\C[\d]$-span of
products of less than $r$ generators with any positioning of brackets and
any orders of multiplication.  

Since the powers of $\d$ can be gathered in the beginning of conformal
monomials (with a probable change in the orders of multiplications), it is
clear that $\bigcup_r C_r=C$.  For a given ordered collection of
generators and a given positioning of brackets, the number of non-zero
monomials is finite because of locality (C1).  Therefore, $\rk C_r$ is
finite.

\df{2.1} The Gelfand-Kirillov dimension of a finitely generated conformal
algebra $C$ is $$\gk
C=\limsup_{r\to\infty}{\log \rk_{\C[\d]} C_r\over \log r}.$$\enddf

Moreover, when $C$ is an associative conformal algebra, every element can
be written as
a sum of $(\dots( \c {f_{j_1}}{n_1}{f_{j_2}}\c ){n_2}{f_{j_3}}\dots\c
){n_{r-1}}{f_{j_r}}$ over $\C[\d]$ with the rewriting process prescribed
by (1.8) and (1.9) not increasing the number of generators involved in the
original presentation of the given element.  Therefore,
$$\eqalign{C_1=&\span_{\C[\d]}(f_1,\dots,f_n),\cr
C_r=&\span_{\C[\d]}(\c gm{f_j}\,|\, g\in C_{r-1}, m\geq 0, 1\leq j\leq
n)+C_{r-1}.}\tag 2.1$$  Hereafter this description of $C_r$'s will be
used.

We will now show that in associative conformal algebras the orders of
multiplications in the monomials used in the presentation (2.1) are
uniformly bounded.  The following lemma is a well-known fact (actually, it
can be deduced directly from the standard proof of Dong's lemma 1.3); the
proof is provided only to demonstrate the employment of rules (1.8) and
(1.9).

\proclaim{Lemma 2.2} Let $N=\max_{i,k} N(f_{j_i},f_{j_k})$. If $(\dots(
\c
{f_{j_1}}{n_1}{f_{j_2}}\c ){n_2}{f_{j_3}}\dots\c ){n_{i-1}}{f_{j_i}}\neq 0$,
then $n_j\leq N$ for all $j$. \endlemma

\pf Denote $(\dots( \c {f_{j_1}}{n_1}{f_{j_2}}\c ){n_2}{f_{j_3}}\dots\c
){n_{i-3}}{f_{j_{i-2}}}$ by $g$. If $n_{i-1}>N$, then $$\eqalign{(\c
g{n_{i-2}}{f_{j_{i-1}}}\c){n_{i-1}}{f_{j_i}}&=\sum_{s\geq n_{i-1}} \binom
{n_{i-2}}{s-n_{i-1}} \c g{n_{i-2}+n_{i-1}-s}(\c {f_{j_{i-1}}}s{f_{j_i}})\cr 
&=\sum_s
\binom {n_{i-2}}{s-n_{i-1}} \c g{n_{i-2}+n_{i-1}-s}0=0.}$$  The statement 
follows by induction.\endpf

One can also speak sometimes of the growth of $C$ meaning the growth of
function $\gamma_C(r)=\rk_{\C[\d]} C_r$.  The Gelfand-Kirillov dimension
of $C$ is finite if and only if $\gamma_C(r)$ is polynomial.
Lemma 2.2 implies that $\rk C_r\leq\sum_{j=1}^n j^rN^{r-1}$. Hence,
just as in the non-conformal case, $\gamma_C(r)$ can not be
superexponential while exponential growth is possible (e.g. in free
conformal algebras, see \cite{Ro} for the definition and explicit
construction of basis).

\rem{2.3} Similar results can be proven for Lie conformal algebras since
for them the presentation (2.1) holds as well because of the Jacobi
identity (1.7).  However the orders of multiplications in the Lie version
of Lemma 2.2 will depend linearly on $r$ (\cite{K1}, \cite{Ro, 1.17}.
\endrem
\ms

The Gelfand-Kirillov dimension is invariant to the change of the
generating set: the new generators are contained in some $C_k$ which in
turn is contained in some $C_l'$ (here $C_1',C_2',\dots$ are the
submodules defined by the new set of generators), thus for the new
$C_r'$'s, $C_r'\subseteq C_{kr}\subseteq C_{lr}'$ and the Gelfand-Kirillov
dimension which measures only the growth of $\gamma_C(r)$ remains the
same.  
Also, the Gelfand-Kirillov dimension of a subalgebra does not exceed the
Gelfand-Kirillov dimension of the full algebra.

Just as in the case of non-finitely generated associative algebras, it is
possible to define a Gelfand-Kirillov dimension of a non-finitely
generated conformal algebra: $$\gk C=\sup_{C'\subset C, C'\text{\ 
finitely generated}}\gk C'.$$

If $C$ is finite as a $\C[\d]$-module, $\gk C=0$.  Conversely, if $C$ is
finitely generated as a conformal algebra and is an infinite
$\C[\d]$-module, it has a non-zero Gelfand-Kirillov dimension. 

We will now relate the Gelfand-Kirillov dimension of $C$ to that of $\coef
C$.  Recall that $\tilde C=C\otimes_{\C}\C[t,t^{-1}]$ and $\coef C=\tilde
C/\tilde\d\tilde C$ where $\tilde\d=\d+\partial/\partial t$.  The map
$(\tilde C,\encirc{0})\to\coef C$ is again denoted by $\phi$.

Notice that by definition, $\tilde\d(ft^k)=\d ft^k+kft^{k-1}$ where $\d
ft^k$ is a shorthand for $(\d f)t^k$.  Therefore, $\phi(\d
ft^k)=-k\phi(ft^{k-1})$. In general, if an element of $\coef C$ is an
image of $\d^ift^k$, one can choose another preimage for it, $f't^{k'}$,
$f'\not\in \d C$. Remark also, that $\dim_{\C} \phi(C')\leq \rk_{\C[\d]}
C'$ for $C'\subseteq C$. 

It is not difficult to see that passing from $C$ to $\tilde C$ increases
$\gk$ by $1$. However, to figure out the effect of factoring out $\tilde
\d\tilde C$ and forgetting the contribution of all multiplications 
but $\encirc{0}$ requires more work.

\proclaim{Theorem 2.4} For an associative conformal
algebra $C$,\newline $\gk \coef C\leq\gk C$+1.\endthm

\pf Notice that even when $C$ is finitely generated, $\coef C$ does not
have to be.  Nonetheless, since we need to prove an upper bound on $\gk
\coef C$, it suffices to demonstrate that such a bound holds for any
finitely generated subalgebra of $\coef C$.

Let $V$ be a generating subspace of $\coef C$.  One can choose
$f_1,\dots,f_n\in C$ and $M^-, M^+\in\Z$ such that
$V\subseteq\span_\C(\phi(f_jt^k)\,|\,1\leq j\leq n, M^-\leq k\leq
M^+).$  We can always assume that $M^-\leq 0$.

Let $$\tilde V=\span_\C(f_jt^k\,|\,1\leq j\leq n, M^-\leq k\leq
M^+)\subseteq\tilde C\tag 2.2.$$ As
$V\subseteq \phi(\tilde V)$, it suffices to prove the upper bound on
growth of $\coef C$ for the subalgebra generated by $\phi(\tilde V)$.
We can also assume that the conformal algebra is smaller, in particular
that it is generated by $f_i$'s.

Consequently, we change notations and take $C$ to be the conformal algebra
generated by $f_1,\dots,f_n$.  Put $V=\phi(\tilde V)$ where $\tilde V$ is
defined as in (2.2).  We shall also use the explicit description of
$C_r$'s given by (2.1).

Put $N=\max_{j_1,j_2} N(f_{j_1},f_{j_2})$ and $M_r^+=rM^+,
M_r^-=r(M^--N)$.  Recall that $M^-\leq 0$.  It is clear that as
functions of $r$, $M_r^+$ increases and $M_r^-$ decreases.

We shall study the growth of $\dim V^r$ via the growth of certain
subspaces of ${\tilde V}^r$ in $(\tilde C,\encirc{0})$.  The inductive
statement claims that
for every element of $V^1+\dots+V^r$ one can choose a preimage in the
subspace ${\tilde V}_r$ of $\tilde C$ spanned by $ht^k$ such that $h\in
C_r$ and $M_r^-\leq k\leq M_r^+$.  (Simultaneously, we will prove that
${\tilde V}^r\subseteq {\tilde V}_r$).

When passing from $r$ to $r+1$, the statement automatically holds for
elements of $V^i,1\leq i\leq r$ as $C_r\subseteq C_{r+1}$ and
$(M_r^-,M_r^+)\subseteq (M_{r+1}^-,M_{r+1}^+)$. Therefore, to prove the
claim one needs only to provide a procedure for choosing the preimage of
an element in $V^{r+1}=V^rV^1$. 
By induction we can consider a larger space, namely
$\phi({\tilde V}_r)V$ and work only with the basis elements of
$\phi({\tilde V}_r)$ and $V$.

Let $gt^{k_1}\in {\tilde V}_r$ where $M_r^-\leq k\leq M_r^+$ and $g\in
C_r$ is a product of $f_j$'s. Consider $f_{j_2}t^{k^2}\in \tilde V$. We
have from (1.5) 
$$\c {gt^{k_1}}0{f_{j_2}t^{k_2}}=\sum_{j\geq 0} \binom
{k_1}j (\c gj{f_{j_2}})t^{k_1+k_2-j}.$$ 
Notice that $j\leq N$ by lemma 2.2.
Therefore $\phi(\c {gt^{k_1}}0{f_{j_2}}t^{k_2})$ lies in
$$\eqalign{&\phi(\span_\C(\c gj{f_{j_2}}t^{k_1+k_2-j}\,|\,M_r^-\leq
k_1\leq M_r^+, M_1^-\leq k_2\leq M_1^+, j\leq N))\subseteq\cr
\subseteq&\phi(\span_\C(ht^k\,|\,h\in C_{r+1},M_r^-+M^--N\leq k\leq
M_r^++M^+)=\cr 
=&\phi(\span_\C(ht^k\,|\,h\in C_{r+1},M_{r+1}^-\leq k\leq
M_{r+1}^+)={\tilde V}_r.}\tag 2.3$$

The immediate consequence is that $\dim V^1+\dots+V^{r+1}$ is bounded by
the dimension of the subspace of ${\tilde V}_{r+1}$ given in (2.3).  It
does not exceed $\dim {\tilde V}_{r+1}=(\rk
C_{r+1})(M_{r+1}^+-M_{r+1}^-)$. We conclude that $$\dim V^1+\dots+V^r\leq
(M_1^+-M_1^-+N')r\cdot\rk C_r$$ and the theorem follows from the
definition of $\gk$.\endpf

\rem{2.5} A similar result was proven for Lie conformal algebras of
growth $0$ in \cite{DK, Lemma 5.2}.\endrem

\ex{2.6} Let $C$ be a conformal algebra which is a torsion
$\C[\d]$-module.  Clearly $\gk C=0$, furthermore, it follows from (C2)  
that any finitely generated subalgebra of $C$ is finite over $\C$. For any
$f\in C$, let $\alpha_i$ be such that $(\Pi_{i=0}^n (\d+\alpha_i))f=0$.  
Put $f_j=(\Pi_{i=j}^n(\d+\alpha_i))f$.  It is clear from (1.6) that the
coefficients of $f_1$ are proportional to either $f_1(0)$ or $f_1(-1)$,
namely for all $n$, $\alpha_0f_1(n)=nf_1(n-1)$.  By induction, the
coefficients of $f$ are linear combinations of $f_j(0)$ and $f_j(-1)$,
$1\leq j\leq n-1$.  Therefore, the coefficient algebra of a finitely
generated subalgebra of $C$ spanned over $\C$ by some $g_i$'s is spanned
over $\C$ by $\{(g_i)_j(0), (g_i)_j(-1)\}$ and is finite over $\C$.  It
follows that $\gk \coef C=0$.  This example shows that the inequality in
Theorem 2.4 is sometimes strict.\endex

\head 3. Unital Conformal Algebras.\endhead

We begin by introducing the analogue of identity in ordinary associative
algebras.

\df{3.1} An element $e\in C$ is called a conformal identity if for every
$f\in C$ $\c e0f=f$  and $N(e,e)=1$.  A conformal algebra containing
such an element is called unital.\enddf

\rem{3.2} Conformal identity is not unique.  Consider, for example, the
current algebra over a unital algebra $A$ that contains a nilpotent
element $r$, $r^2=0$.  Then both $f_{\bold 1}$ and $f_{{\bold
1}+r}=f_{\bold 1}-\d f_r$ are conformal identities in $\cur(A)$ (here
${\bold 1}$ is the identity in $A$).  An obvious generalization of this
example will be used in the proof of Theorem 4.2.\endrem
\ms

As we will demonstrate immediately, the Weyl algebra fits the definition.
Furthemore, in a sense it is the only such simple algebra of growth $1$
(see Theorem 4.6 for details). 

It is one of the nice properties of Lie conformal algebras that every
torsion element is central (\cite{DK, Prop. 3.1}).  In general, this is
not true for associative algebras; however, for unital algebras an even
stronger result holds.

\proclaim{Lemma 3.3} A unital conformal algebra $C$ is always
torsion-free.\endlemma

\pf The statement can be deduced from \cite{DK}; however, an easier
proof is possible.  Let $f\in C$ be a torsion element. Since $\C$ is
algebraically closed, we can always assume without loss of generality that
$(\d+\alpha)f=0$ for some $\alpha\in\C$.  Clearly for all $n>0, g\in C$,
$(\d+\alpha)(\c gnf) =-n\c g{n-1}f$.  If $\c g0f\neq 0$, $\c gnf\neq 0$
for all $n\geq 0$ contradicting locality of $g$ and $f$.  Thus, $C$ is the
left annihilator of $f$ which contradicts the definition of $e$.\endpf

It is an open question whether every torsion-free conformal algebra can be
embedded into a unital conformal algebra.  This question can be
reformulated in terms of the coefficient algebras as well (see Corollary
to Proposition 3.5).  We remark that such an embedding automatically
exists if the conformal algebra possesses a finite faithful representation
as it implies that the algebra can be embedded into $\text{CEnd}_N$
\cite{K1, ch.2}.\ms

We now turn to the description of unital associative conformal algebras.
The following fact is well-known (see \cite{Ro} for a fuller
exposition). The proof is elementary; however, some of the explicit
calculations it contains will be needed below.

\proclaim{Lemma 3.4} Let $A$ be an associative algebra, $\delta$ a
locally nilpotent derivation of $A$.  Then the set $F=\{f_a=\sum at^n
z^{-n-1}, a\in A\}$ of series over a localized Ore extension $A[t,
t^{-1};\delta]$ of $A$ is a conformal algebra with the operations $\c
{f_a}m{f_b}=f_{(-1)^ma\delta^m(b)}.$\endlemma

\pf Let $a, b$ be two elements of $A$.  Consider the series $f_a=\sum at^n
z^{-n-1}$ and $f_b=\sum bt^n z^{-n-1}$.  The $k$-th coefficient of their
$m$-th product is $\sum_s (-1)^s \binom ms at^{m-s}bt^{k+s}$. Since
$t^nb=\sum_i (-1)^i \binom ni\delta^i(b)t^{n-i}$, we have
$$(\c {f_a}m{f_b})(k)=\sum_{s,i} (-1)^{s+i}\binom {m-s}i\binom ms
a\delta^i(b)t^{m+k-i}.$$
If $i\neq m$, $\sum_s (-1)^{s+i}\binom {m-s}i\binom ms=0$,
thus $(\c {f_a}m{f_b})(k)=(-1)^ma\delta^m(b)t^k$. Hence,
$$\c {f_a}m{f_b}=f_{(-1)^ma\delta^m(b)}.\tag 3.1$$
It follows immediately that $F$ is closed with respect to taking conformal
products. Locality follows from the nilpotency of $\delta$.\endpf

Conformal algebras described above are called {\it differential algebras}.

If one takes $A=C[x]$ and $\delta=\d/\d_x$, the resulting conformal
algebra built as above is precisely the Weyl algebra ${\hbox {\got
W}}$.

Notice that for the differential algebras constructed in Lemma 3.4, $\gk
C=\gk \coef C-1$. Indeed, the only osbtruction to the equality in the
proof of Theorem 2.4 are non-zero products of generators with a zero
zeroeth coefficient.  However, if a formal distribution $f_a$ has
$f_a(0)=0$, then it is zero itself. 

We will now establish the converse of Lemma 3.4.

It should be noted first that for a conformal identity $e$, $e(0)$ is a
left but not necessarily right identity of $\coef C$.  Clearly, if the
left annihilator of an associative algebra is $0$, the left identity is
also a right identity.  We need to introduce a similar condition for
conformal algebras.  We call the set $$L(C)=\{f\in C\,|\,\c fng=0 \text{\
for all\ } g\in C,n\in\Z\}$$ the left annihilator of $C$.  Notice that due
to (1.8) and (1.9) $L(C)$ is a nilpotent ideal of $C$.

\proclaim{Proposition 3.5} Let $C$ be a unital conformal algebra with the
conformal identity $e$.  Let $C$ have a trivial left annihilator. Then
$\coef C=A[t,t^{-1};\delta]$ where $\gk A=\gk C$ and $\delta$ is a locally
nilpotent derivation of $A$.  If $C$ is finitely generated, then so is
$A$.\endprop

\pf Notice first that $(\c e1e)(0)=e(1)e(0)-e(0)e(1)=0$.  Thus $e(0)$ is
an identity in the subalgebra of $\coef C$ generated by $e(0)$ and $e(1)$.
Later we will show that $e(0)$ is an identity in $\coef C$.

Consider, as before, the algebra $\coef C$ as an image of
$$\phi:C[t,t^{-1}]\to C[t,t^{-1}]/(\d+\d_t) C[t,t^{-1}].$$ 
Let $A$ be the
image of $C$ in $\coef C$.  For the conformal identity $e$, (1.5) implies
$$et^{n+1}=\c {et^n}0{et}-\sum_{j=1}^n\binom nj(\c eje)t^{n+1-j}=\c
{et^n}0{et}.$$ By induction, $et^n=(et)^n$, $n>0$. Since $\c
{et^n}0{et^{-n}}=\c e0e$, we have
$$e=\sum_{n\in \Z} t^nz^{-n-1}.\tag 3.2$$
(With a slight abuse of notations, we denote $\phi(et)$ by $t$.) In these
notations $t^0=e(0).$

It will follow that $\coef C$ is generated by the algebra of $0$-th
coefficients $A=\phi(C)$, $t$ and $t^{-1}$ if one can choose a set of
generators $\{g_i\}$ of $C$ such that $\c {g_i}0e=g_i$.  Notice that for
such elements $g(0)e(0)=g(0)$, thus it will follow that $e(0)$ is indeed
the identity of $\coef C$. In general, in terms of coefficients
$g_i(n)=g_i(0)t^n$ for all $n$. We will now show how every element of $C$
can be expressed via the elements of this particular type.

Assume at first that $e(0)$ is an identity in $\coef C$.

Consider now an arbitrary element $f=\sum f(n)z^{-n-1}$. Put
$\widehat{f(n)}= f(n)t^{-n}$.  From the locality condition (1.4) in Lemma
1.1, we see that
$$\sum_{j=0}^{N(f,e)+1} (-1)^j\binom {N(f,e)+1}j \widehat{f(n-j)}=0\ \ \
\text{for all }n>N(f,e),\tag 3.3$$
as $f(n-j)e(j)=\widehat{f(n-j)}t^{n-j}\cdot t^j$.  It follows that for any
$n$, $\widehat{f(n)}$ is determined by a polynomial of degree $N(f,e)$
over $\coef C$, namely
$\widehat{f(n)}=f_0+f_1n+\dots+f_{N(f,e)}n^{N(f,e)}$. Indeed, if
$N(f,e)=0$, then $\widehat{f(n)}=\widehat{f(n-1)}$ and the result is
clear.  In the case of arbitrary $N(f,e)$, the locality condition (3.3)
can be rewritten as
$$\sum_{j=0}^{N(f,e)}(-1)^j \binom {N(f,e)}j
(\widehat{f(n-j)}-\widehat{f(n-j+1}))=0\tag 3.4$$
follows by induction.

Thus, we can rewrite $f$: 
$$f=\sum (f_0+f_1n+\dots+f_{N(f,e)}n^{N(f,e)})t^n z^{-n-1}.$$  

If $N(f,e)=0$, the element $f=\sum f_0t^n z^{-n-1}$ already has the form
we need.  In the general case, consider the
product $\c fNe$, $N=N(f,e)$. Its coefficients are
$$(\c fNe)(n)=\sum (-1)^j \binom Nj 
(f_0+f_1(N-j)\dots+f_N(N-j)^N)t^{N-j}t^{n+j}.\tag 3.5$$
In the expression above the coefficient at $f_kt^{N+n}$ is 
$$\sum (-1)^j \binom Nj (N-j)^k=\Biggl\{\eqalign{&\ 0,\ \ \ \ k<N\cr
&\ N!,\ \ k=N}.$$
Therefore, $\c fNe=\sum (N!f_Nt^N)t^n z^{-n-1}$ and has the form we need.
It immediately follows that
$\d^N(\c fNe)=\sum (-1)^Nn(n-1)\cdots (n-N)N!f_Nt^Nz^{-n-1}$
and the $n$-th coefficient of $$f-(-1)^N{1\over N!}\d^N(\c fNe)$$
is a polynomial in $n$ of degree $N-1$ with constant coefficient.  It is easy
to see that its degree of locality with $e$ is also at most $N-1$. Hence,
by repeating the process described above we can express $f$ as a linear
combiantion of monomials of elements in the desired 
form.  By taking as $f$ each generator of $C$ we obtain a set of
generators $\{g_i\}$ such that for all $i$ and $n$, $g_i(n)=g_i(0)t^n$.
Moreover, if $C$ is finitely generated, this set is also finite.

We now return to the general case (i.e. $e(0)$ is not necessarily an
identity).  Consider a map $\Xi:C\to\coef C[[z,z^{-1}]]$ where
$\Xi(f)=\sum f(n)e(0)z^{-n-1}$.  Denote the image of $\Xi$ by $C'$.
Clearly $\Xi:C\to C'$ is a linear map.  Moreover, since $e(0)$ is a left
identity, $f(k)g(m)e(0)=(f(k)e(0)g(m)e(0)$ and (1.2) implies $\Xi(\c
fng)=\c {\Xi(f)}n{\Xi(g)}$. Hence $C'$ is a unital conformal algebra of
formal distributions over $(\coef C)e(0)$. We do not claim here that
$\coef C'=(\coef C)e(0)$.  However, the above construction of a particular
generating set $\{g_i\}$ of a unital conformal algebra does not use the
universality of a coefficient algebra. Hence, $C'$ has such a system of
generators: $\{g_i'=\sum g_i'(0)t^nz^{-n-1}\}$.

Notice now that $\Xi:C\to C'$ is a bijection.  Indeed, if $\Xi(f)=0$, it
follows that since $f(k)e(0)=0$ for all $k$, $f(k)g(m)=0$ for all $g\in
C, m,k\in\Z$.  Thus $\c fng=0$ and $f\in L(C)$, a contradiction.

Put $g_i=\Xi^{-1}(g_i')$, $g_i'$ as above.  Notice that $g_i$'s are
generators of $C$. Then $g_i(n)e(0)=g_i'(0)t^n=g_i(0)e(0)t^n=g_i(0)t^n=(\c
{g_i}0e)(n)$.  Thus $\Xi(g_i)=g_i'=\c {g_i}0e\in C$ and $C'$ is a
subalgebra of $C$. Moreover, for every $f\in C, n\in\Z_{\geq 0}$, (1.2)
implies that $\c {(g_i-\Xi(g_i))}nf=0$.  Hence $\Xi(g_i)=g_i$ and $C'=C$.

We conclude that $C$ possesses a system of generators $\{g_i\}$ such that
$g(n)=g(0)t^n$.  In particular this implies that $\coef C=A[t,t^{-1}]$
with identity $e(0)$.

Define now a derivation $\tilde\delta(f)=\c f0{et}-\c {et}0f$ in $(\tilde
C, \encirc{0})$. Then, on generators
$$\tilde\delta(f)=\c f0{et}- (\c e0{ft}+\c e1f)=-\c e1f.\tag 3.6$$
Next, we introduce a derivation
$\delta(\phi(f))=\phi(\tilde\delta(f)$ on $A$.
In this setup, $\coef
C$ is an Ore extension of $A$. Moreover, $$\eqalign{
\tilde\delta^n(f)=&\c e1{(\c
e1{\tilde\delta^{n-2}f})}=\cr =&\c {(\c e0e)}2{\tilde\delta^{n-2}f}-\c
{(\c e1e)}1{\tilde\delta^{n-2}f}=\c e2{\tilde\delta^{n-2}(f)}}$$ and by
induction, $\tilde\delta^nf=(-1)^n\c enf$. Therefore, $\delta$ is locally
nilpotent.\ms

It follows from the calculation in Lemma 3.4 that if $C$ is finitely
generated, then so are $A$ and $\coef C$.  In particular, the generators
of $\coef C$ are the zeroeth coefficients of the generators of $C$, their
$\delta$-derivatives, and $t$.  For unital algebras, $\phi(f_a)=0$ if and
only if $f_a=0$.  Thus, if in the proof of Theorem 2.4 one chooses the
generating subspace of $\coef C$ as described above, all inequalities turn
into equalities, hence $\gk\coef C=\gk C+1$.

By \cite{KL, 3.5 and 4.9} $\gk A=\gk A[t,t^{-1};\delta]-1=\gk C$.  \endpf

\proclaim{Corollary} An associative conformal algebra $C$, such that
$L(C)=0$, can be embedded into a unital algebra if and only if it its
coefficient algebra can be embedded into a localized Ore extension defined
by a locally nilpotent derivative.\endcor

Nonisomorphic unital conformal algebras $C_1$ and $C_2$ yield
nonisomorphic pairs ($A_1,\delta_1$) and ($A_2,\delta_2$).  Indeed, even
if $A_1\cong A_2$, then $\delta_1\neq\delta_2$, otherwise by Lemma 3.4
and, in particular, (3.1), there exists a canonical isomorphism between
$C_1$ and $C_2$.  Clearly non-isomorphic $A_1$ and $A_2$ give rise to
non-isomorphic $C_1$ and $C_2$; however, pairs $(A,\delta_1)$,
$(A,\delta_2)$ with $\delta_1\neq\delta_2$ may lead to the same conformal
algebra.  Consider, for example, $A=\Mat_n(\C[x])$ with $\delta_1=\d/\d
x$, and $\delta_2=\d/\d x+\ad a$ where $a\in A$ is a nilpotent matrix (the
construction of an isomorphism between $(A,\delta_1)$ and $(A,\delta_2)$
follows the proof of Theorem 4.6).

\head 4. Unital Conformal Algebras of Small Growth.\endhead

The most immediate implication of Proposition 3.6 is that any
classification of unital conformal algebras is much easier than the
general case because one needs only to describe the associative algebra
$A$ and its derivation $\delta$. Here we will do just that for conformal
algebras of growth $0$ and $1$ over $\C$ (Actually, any algebraically
closed field of characteristic 0). For the rest of this chapter, $C$
always denotes a unital associative conformal algebra with the zero left
annihilator, and $A$ and $\delta$ are the resulting associative algebra
and its locally nilpotent derivation.

\proclaim{Lemma 4.1} $C$ is simple if and only if $A$ contains no
$\delta$-stable ideals.\endlemma

\pf If $I$ is a $\delta$-stable ideal of $A$, then $I[t,t^{-1}]$ is an
ideal of $\coef C$ and $J=\{f_a\,|\,a\in I\}$ an ideal of $C$ by
calculations in Lemma 3.4.

Conversely, let let $J$ be a an ideal of $C$ and put $I=\phi(J)$, the set
of all $0$-th coefficients of $J$. Clearly $a\in I$ if and only if $f_a\in
J$. Since for any $a\in I$, $b\in A$, $f_{ab}=\c {f_a}0{f_b}\in J$, $I$ is
an ideal of $A$.  Moreover, since $f_{\delta^n(a)}=(-1)^n(\c en{f_a})\in
J$, $I$ is a $\delta$-stable ideal of $A$.

It remains to show that $I$ is non-empty.  Since $C$ is a coefficient
algebra, for any $f\in C$, there exists a collection $\{a_i\}\subset A$
such that $f=\sum_{i,k}c_{ik}\d^kf_{a_i}$.  Therefore, for every $m$,
$$\eqalign{f(m)=&\sum_{i,k} c_{ik}(\d^kf_{a_i})(m)
=\sum_{i,k}(-1)^m{m!\over (m-k)!} c_{ik}f_{a_i}(m-k)=\cr
=&\sum_k \Big(\sum_i c_{ik}a_i\Big)\Big((-1)^m {m!\over
(m-k)!}t^{m-k}\Big)}\tag 4.1$$
(if $m\geq 0$, we sum only over $k\leq m$).  It follows that if $f(m)=0$
for all $m\geq 0$, then $\sum_i c_{ik}a_i=0$ for all $k$ implying
$f(m)=0$ for $m<0$.

Let now $f$ be an element of $J$, and $m\geq0$ minimal such that
$f(m)\neq0$.  Then $(\c fme)(0)\neq 0$.
\endpf

\proclaim{Corollary} Let $C$ be an associtive current algebra $\cur(A)$
(not necessarily unital) where $A$ contains more than $1$ element. Then
$C$ is simple iff $A$ is simple.\endcor

\pf The proof is the same as the proof of the Lemma; it is only required
that $A$ contain $a$ such that $\c fm{f_a}\neq 0$, $f$ and $m$ as above.
If $f(m)A=0$, then there exists $b\in A$, $bA=0$. It follows that $A$ is
not simple.\endpf

\proclaim{Theorem 4.2} If $C$ is a simple unital associative conformal
algebra and a finite $\C[\d]$-module, then it is a current algebra over
the ring of matrices over $\C$. \endthm

\pf By Proposition 3.5 $\gk A=0$ and Lemma 4.1 $A$ is $\delta$-simple. It
follows immediately from \cite{Bl} that $A=\Mat_n(\C)$.  It is well-known
that $\delta$ is necessarily an inner derivation $\ad a$(see, e.g.
\cite{H}). It is easy to show that there exists a nilpotent $r\in A$ such
that $\delta=\ad r$ ($r$ is the nilpotent part of $a$).

Clearly, the ring of skew-polynomials $A[t, t^{-1};\delta]$ is isomorphic
to the ring of polynomials $A[s,s^{-1}]$ (send $t-r$
to $s$). To prove that the corresponding differential algebras are
isomorphic
as well, it suffices to show that the series $e'=\sum (t-r)^nz^{-n-1}$
belongs to the differential algebra over $A[t,t^{-1};\delta]$.  Indeed, in
this case the conformal algebra generated by $e'$ and $\c {\big(\sum
at^nz^{-n-1}\big)}0e'$ is canonically isomorphic to $\cur(A)$.

Let $m$ be the degree of nilpotency of $r$. Since $\delta(r)=0$, $t$ and
$r$ commute, therefore, $$e'=\sum_{k=0}^{m-1} {1\over k!} \d^k \Big( \sum
r^kt^nz^{-n-1}\Big) \in \text{diff. algebra over\ }A[t,t^{-1};\delta].$$
\endpf

The classification of simple associative algebras of finite type (i.e.  
those finite as $\C[\d]$-modules) \cite{K2, 4.4} implies that all simple
algebras of finite type are unital.  However, so far there exists no
direct way to demonstrate that a simple algebra of finite type can be
adjoined a conformal identity.\ms

The next task is to classify simple unital conformal algebras of growth
$1$. Of great help is the classification of finitely generated associative
algebras of growth $1$ developed in \cite{SW} and \cite{SSW}.  We will use
the following result of the former:

\proclaim{Theorem \cite{SW}} A finitely generated prime algebra $A$, $\gk
A= 1$ is necessarily a finite module over its center which is also
finitely generated. \endthm

\rem{4.3} The immediate corollary of this classification is an otherwise
non-obvious fact that there exist no finitely generated simple associative
current conformal algebras of linear growth.  Indeed, let $C=\cur(A)$ be
such an algebra.  $A$ is simple as well by corollary to Lemma 4.1 and $\gk
A=1$.  Hence, its center is a finitely generated simple algebra of $\gk
1$, i.e. a field of transcedental degree $1$.  Such can not be finitely
generated over $\C$ as an algebra.\endrem
\ms

Assume now that $C$ (hence $A$) is finitely generated. The follwong
lemma is well-known.

\proclaim{Lemma 4.4} \cite{Po} Let $A$ be a differentially simple algebra
as above. Then $A$ is prime. \endlemma

We begin the classification of unital conformal algebras of growth $1$
(i.e. the classification of pairs ($A$,$\delta$) up to the isomorphism of
resulting conformal algebras) with the commutative case. 

\proclaim{Lemma 4.5} If $A$ is commutative (not necessarily finitely
generated), then either $\delta\equiv 0$ and $A$ is a field of
transdendental degree one or $A=\C[x]$ and $\delta=\partial/\partial x$
for some $x\in A$. \endlemma

\pf Consider an element whose derivative is $0$.  It generates a
$\delta$-stable ideal of $A$, hence it is invertible and the set of all
such elements forms a subfield of $A$.

It follows that the set of algebraic elements of $A$ whose derivative is
$0$ forms an algebraic subfield of $A$, i.e. is $\C$.  Let $x\in A$ be
transcendental and $n$ minimal such that $\delta^n(x)=0$.  Two cases can
occur: either $\delta^{n-1}(x)\in\C$ (we can always assume that it is $1$)
or $\delta^{n-1}(x)$ is transcendental. Hence, the set of all
transcendental elements of $C$ splits into two, each corresponding to one
of the cases above.  We will show that these cases do not occur
simultaneously (i.e. either one of the subsets is necessarily empty). Pick
an arbitrary element from each set, respectively $x$ and $y$. Without loss
of generality we may assume that $\delta(x)=1,\delta(y)=0$.  These two
elements are algebraically dependent since $\gk A=1$ but in every
statement of dependence of $y$ over $\C[x]$ the degree in $x$ can be
lowered by the application of $\delta$ making $y$ algebraic over $\C$.  
Thus these cases do not occur simultaneously.

If we assume that no transcendental element has a non-zero derivative,
i.e. $\delta\equiv 0$, $A$ is a field of transcendental degree $1$ over
$\C$.  Moreover, if at least one transcendental element has a zero
derivative, we can apply the argument from the last paragraph: if there
exists an element with a non-zero derivative, then there exist $x$ and $y$
be such that $\delta(x)=y$, $\delta(y)=0$.  But $y$ is invertible and
$\delta(y^{-1}x)=1$. We conclude that if one transcendental element has a
zero derivative, then every transcendental element does and $A$ is a
field.

The remaining case follows from a lemma in \cite{W}; however, in the
current framework one can provide an elemetary proof.  If there exists $x$
such that $\delta(x)=1$, consider an arbitrary $y\in A$ with
$\delta(y)=x$. We have $\delta(y-x^2/2)=0$ which by the previous
discussion means that $y=x^2+a, a\in\C$ (i.e. $y$ is an ``antiderivative''
of $x$).  For an arbitrary $y\in A$ we conclude by induction on $n$ such
that $\delta^n(y)=0$, that $y\in\C[x]$.\endpf

We already know that for a simple unital conformal algebra $C$, $\coef
C=A[t,t^{-1};\delta]$ where $\delta$ is a locally nilpotent derivative and
$A$ is a prime algebra of $\gk$ $1$.  It follows from \cite{SW} that $A$
is a finite module over its center $Z(A)$ which is also of $\gk$ $1$.
Notice that $\delta$ can always be restricted to $Z(A)$ which is
$\delta$-simple as well. Indeed, assume the contrary, i.e. the existence
of $a\in Z(A)$ that, together with its derivatives, generates a
$\delta$-stable ideal of $Z(A)$.  Since $\delta$ is locally nilpotent,
without loss of generality we may assume that $\delta(a)=0$ and that
$\langle a\rangle$ is an ideal of $Z(A)$.  But by $\delta$-simplicity of
$A$, there exists $b\in A$, $ab=1$, thus for any $x\in A$,
$xb-bx=ba(xb-bx)=b(x-x)=0$.  Hence, $a$ is invertible in $Z(A)$, a
contradiction.  Thus we can apply lemma 4.5 to the classification of the
general case.

\proclaim{Theorem 4.6} Let $C$ be a finitely generated simple unital
associative conformal algebra of linear growth. Then $C$ is isomorphic to
the differential algebra of formal distributions over
$\Mat_n(W_t)$.\endthm

\pf From the Proposition 3.5 and Lemma 4.4, $\coef C=A[t,t^{-1};\delta]$,
$Z(A)=\C[x]$, $\delta |_{\C[x]}=\d/\d x$.

Consider subalgebra $A_0={\ker}\delta$.  We begin by demonstrating that
$A_0$ generates $A$ over $Z(A)$. More precisely, every $a\in A$ is of the
form $\sum_1^n {x^i\over i!}a_i$, $a_i\in A_0$ where $n$ is such that
$\delta^n(a)=0, \delta^{n-1}(a)\neq 0$. Indeed, with the induction
assumption $\delta(a)=\sum_1^{n-1} {x^i\over i!} b_i$, $b_i\in A_0$,
consider $a_0=a-\sum_1^{n-1}{x^{i+1}\over (i+1)!}b_i$. Since
$\delta(a_0)=0$, we see that $a$ is also a polynomial in $x$ over $A_0$. 
Moreover, this polynomial form is unique for any $a\in A$, i.e. if
$\sum_1^n x^ia_i=0$, applying a necessary number of derivatives
shows that the coefficient at the highest power is $0$.

Fix a subset $\{a_i\}$ of $A_0$ that generates $A$ over $Z(A)$. Any
product of elements of $A_0$ is a linear combination $\sum
p_i(x)a_i$ over $Z(A)$ with the derivative $\sum (\d p_i(x)/\d x)a_i=0$. 
This implies $A_0=\span_{\C}(a_i)$ is finite dimensional. 

Clearly any ideal of $A_0$ can be lifted to $A$, thus $A_0$ is prime as
well and therefore simple over $\C$ (\cite{R, 2.1.15}). 
Hence $A=\Mat_n(\C[x])$ and $\delta=\d/\d x$. It follows that
$A[t,t^{-1};\delta]$ is a matrix ring over $\C[x][t,t^{-1};\d/\d
x]=W_t$.\endpf

\enddocument